\documentclass[leqno,12pt]{amsart}
\usepackage{amsmath}
\usepackage{amsfonts}
\usepackage{amssymb}
\usepackage{amsthm}
\usepackage[usenames,dvipsnames]{color}
\usepackage{hyperref}
\newtheorem{thm}{Theorem}
\newtheorem{cor}[thm]{Corollary}
\newtheorem{lem}[thm]{Lemma}
\newtheorem{prop}[thm]{Proposition}
\theoremstyle{definition}
\newtheorem{defn}[thm]{Definition}
\theoremstyle{remark}

\theoremstyle{remark}
\newtheorem{eg}[thm]{Example}
\newenvironment{acknowledge}{\noindent\textbf{Acknowledgments.}}{}
\numberwithin{equation}{section}
\numberwithin{thm}{section}
\long\def\blankfootnotetext#1{\begingroup\def\thefootnote{\fnsymbol{footnote}}\footnotetext{#1}\endgroup}
\newcommand{\Proj}{\mathbb{P}}
\newcommand{\Z}{\mathbb{Z}}
\newcommand{\Q}{\mathbb{Q}}
\newcommand{\R}{\mathbb{R}}
\newcommand{\set}[1]{\left\{{#1}\right\}}
\newcommand{\abs}[1]{\left\vert{#1}\right\vert}

\newcommand{\sconv}[1]{\mathrm{conv}\!\set{#1}}

\newcommand{\scone}[1]{\mathrm{cone}\!\set{#1}}
\newcommand{\Hom}[1]{\mathrm{Hom}({#1})}
\newcommand{\modb}[1]{\mathrm{(mod}\,#1\mathrm{)}}
\renewcommand{\thefootnote}{\fnsymbol{footnote}}
\begin{document}
\title{Bounds on Fake Weighted Projective Space}
\author{Alexander M. Kasprzyk}
\address{Department of Mathematics and Statistics\\University of New Brunswick\\Fredericton NB\\E3B 5A3\\Canada.}
\email{kasprzyk@unb.ca}
\maketitle
\blankfootnotetext{2000 \textit{Mathematics Subject Classification}. Primary 14M25; Secondary 14J45, 52B20.}
\blankfootnotetext{\textit{Key words and phrases}. Weighted projective space, canonical, terminal.}
\begin{abstract}
A fake weighted projective space $X$ is a $\Q$-factorial toric variety with Picard number one. As with weighted projective space, $X$ comes equipped with a set of weights $(\lambda_0,\ldots,\lambda_n)$. We see how the singularities of $\Proj(\lambda_0,\ldots,\lambda_n)$ influence the singularities of $X$, and how the weights bound the number of possible fake weighted projective spaces for a fixed dimension. Finally, we present an upper bound on the ratios $\lambda_j/\sum\lambda_i$ if we wish $X$ to have only terminal (or canonical) singularities.
\end{abstract}
\section{Introduction}\label{sec:fake_wps_definitions}
Let $N\cong\Z^n$ be an $n$-dimensional lattice, and $N_\R:=N\otimes_\Z\R$. The dual lattice $M:=\Hom{N,\Z}\cong\Z^n$ is often referred to as the \emph{monomial lattice}. Let $\set{\rho_0,\rho_1,\ldots,\rho_n}\subset N$ be a set of primitive lattice points such that $N_\R=\sum_{i=0}^n\R_{\geq 0}\rho_i$. There exist $\lambda_0,\lambda_1,\ldots,\lambda_n\in\Z_{>0}$, with $\gcd\!\set{\lambda_0,\lambda_1,\ldots,\lambda_n}=1$, unique up to order, such that:
$$\lambda_0\rho_0+\lambda_1\rho_1+\ldots+\lambda_n\rho_n=0.$$
Define the $n$-dimensional cones:
$$\sigma_i:=\scone{\rho_0,\rho_1,\ldots,\widehat{\rho}_i,\ldots,\rho_n},\qquad\text{for }i=0,1,\ldots,n,$$
where $\widehat{\rho}_i$ indicates that the point $\rho_i$ is omitted. The $\sigma_i$ generate a complete $n$-dimensional fan $\Delta$.
\begin{defn}
The projective toric variety associated with the fan $\Delta$ is called a \emph{fake weighted projective space with weights $(\lambda_0,\lambda_1,\ldots,\lambda_n)$}.
\end{defn}
An immediate consequence of this definition is that fake weighted projective spaces are $\Q$-factorial toric varieties with Picard number one. Of course, the collection of weighted projective spaces is a sub-collection of fake weighted projective spaces. Naturally, there exist fake weighted projective spaces which are not weighted projective spaces.
\begin{eg}
Consider the cubic surface $(W^3=XYZ)\subset\Proj^3$. This has three $A_2$ singularities, and can be realised as $\Proj^2/(\Z/3)$, where the $(\Z/3)$-action is given by:
$$\varepsilon:x_i\mapsto\varepsilon^ix_i\qquad\text{for }(x_1,x_2,x_3)\in\Proj^2,$$
where $\varepsilon$ is a third root of unity. The corresponding fan has rays $\rho_0=(2,-1)$, $\rho_1=(-1,2)$, and $\rho_2=(-1,-1)$; it is a fake weighted projective surface with weights $(1,1,1)$.
\end{eg}
\begin{eg}
Consider the three-dimensional toric variety $X$ generated by the fan with rays $\rho_0=(1,0,0)$, $\rho_1=(0,1,0)$, $\rho_2=(1,-3,5)$, and $\rho_3=(-2,2,-5)$. This is a fake weighted projective space with weights $(1,1,1,1)$, not isomorphic to $\Proj^3$. In fact $(-K)^3=64/5$, and $X$ has four terminal singularities of type $\frac{1}{5}(1,2,3)$.
\end{eg}
The second example above has appeared on several occasions in the literature (\cite[pg.~178]{BB92},~\cite[Remark 14.2.3]{Mat02},~\cite[pg.~189]{Snow05}, and~\cite[Table~4]{Kas03}). An interesting construction can be found in~\cite[\S 4.15]{Reid85}:
\begin{eg}
Let $M\subset\Z^4$ be the three-dimensional affine lattice defined by:
$$M:=\set{(m_1,m_2,m_3,m_4)\in\Z^4\left|\sum_{i=1}^4 m_i=5\text{ and }\sum_{i=1}^4im_i\equiv 0\ \modb{5}\right.}.$$
Let $\Sigma\subset M_\R$ be the simplex whose four vertices are given by the points $(5,0,0,0),\ldots,$ $(0,0,0,5)$ (i.e. the points corresponding to the monomials $x_i^5$, $i=1,2,3,4$). The toric variety constructed from $\Sigma$ is $\Proj^3/(\Z/5)$, where the $(\Z/5)$-action is given by:
$$\varepsilon:x_i\mapsto\varepsilon^ix_i\qquad\text{for }(x_1,x_2,x_3,x_4)\in\Proj^3,$$
where $\varepsilon$ is a fifth root of unity.
\end{eg}
Fake weighted projective spaces occur naturally in toric Mori theory. The following result is adapted from~\cite[(2.6)]{Reid83T}. (The original statement claimed that all the fibres of $\varphi_R\!\mid_A$ were weighted projective spaces. This oversight has been noted -- and corrected -- in, amongst other places,~\cite[Remark~14.2.4]{Mat02},~\cite[\S 1]{Fuj03}, and~\cite{Krych02}.)
\begin{prop}\label{prop:fwps_fibres_Mori_theory}
Let $X$ be a projective toric variety whose associated fan $\Delta$ is simplicial (i.e. $X$ is $\Q$-factorial). If $R$ is an extremal ray of $N\!E(X)$ (the cone of effective one-cycles) then there exists a toric morphism $\varphi_R:X\rightarrow Y$ with connected fibres, which is an elementary contraction in the sense of Mori theory: $\varphi_{R^*}\mathcal{O}_X=\mathcal{O}_Y$, and for a curve $C\subset X$, $\varphi_RC$ is a point in $Y$ if and only if $[C]\in R$.

Let
$$
\begin{array}{r@{\hspace{2mm}}c@{\hspace{2mm}}l}
A&\longrightarrow&B\\
\cap&&\cap\\
\varphi_R:X&\longrightarrow&Y
\end{array}
$$
be the loci on which $\varphi_R$ is not an isomorphism. Then $\varphi_R\!\mid_A:A\rightarrow B$ is a flat morphism, all of whose fibres are fake weighted projective spaces of dimension $\dim A-\dim B$.
\end{prop}

We shall investigate the relation between a fake weighted projective space $X$ with weights $(\lambda_0,\ldots,\lambda_n)$ and $\Proj(\lambda_0,\ldots,\lambda_n)$. In particular, we shall see that the `niceness' of the singularities of $X$ is restricted by the singularities of $\Proj(\lambda_0,\ldots,\lambda_n)$ (Corollaries~\ref{cor:map_fake_to_wps} and~\ref{cor:map_Gorenstein_to_wps}). We shall also introduce a measure of `how much' $X$ differs from the \emph{bona fide} weighted projective space, and establish an upper bound on this measure (Theorem~\ref{thm:bound_fake_wps_mult} and Corollary~\ref{cor:bound_canon_fake_wps_mult}). Finally, in Theorem~\ref{thm:barycentric_bounds}, we present an upper bound on the weights if we wish $X$ to have at worst terminal (or canonical) singularities.

\bigskip\begin{acknowledge}
The author would like to express his gratitude to Dr.~G.~K.~Sankaran for his invaluable explanations and advice, and to Professor Victor Batyrev for his encouragement and observations. The author would also like to thank Weronika Buczy\'nska for introducing him to the name ``fake weighted projective space'' in her talk \emph{Fake weighted projective spaces and Mori theorem for orbifolds}, presented at the Calf Seminar, December 2005.

A significant portion of this work was funded by an Engineering and Physical Sciences Research Council (EPSRC) studentship, and forms part of the author's PhD thesis~\cite{KasPhD}. The author is currently funded by an ACEnet Postdoctoral Research Fellowship.
\end{acknowledge}
\section{Fake Weighted Projective Space and Weighted Projective Space}
We consider what can be said about fake weighted projective space in terms of the corresponding weighted projective space. In particular, we shall see how the singularities of one are dictated by the other, and how weighted projective space provides a bound when searching for fake weighted projective spaces.

We shall rely on the following result, which allows us to distinguish between fake and genuine weighted projective space:
\begin{prop}[\protect{\cite[Proposition~2]{BB92}}]\label{prop:generating_fan_unique}
For any weights $(\lambda_0,\lambda_1,\ldots,\lambda_n)$ such that $\gcd\!\set{\lambda_0,\lambda_1,\ldots,\lambda_n}=1$,  let $\rho_0,\rho_1,\ldots,\rho_n\in N$ be the primitive generators for the fan of $\Proj(\lambda_0,\lambda_1,\ldots,\lambda_n)$. Then:
\begin{itemize}
\item[(i)] $\lambda_0\rho_0+\lambda_1\rho_1+\ldots+\lambda_n\rho_n=0$;
\item[(ii)] The $\rho_i$ generate the lattice $N$.
\end{itemize}

Furthermore, if $\rho'_0,\rho'_1,\ldots,\rho'_n$ is any set of primitive lattice elements satisfying (i) and (ii) then there exists a transformation in $GL(n,\Z)$ sending $\rho_i$ to $\rho'_i$ for $i=0,1,\ldots,n$.
\end{prop}
Note that two complete toric varieties are isomorphic as abstract varieties if and only if they are isomorphic as toric varieties~(\cite{Dem70}). This is a consequence of fact that the automorphism group is a linear algebraic group with maximal torus; Borel's Theorem tells us that in such a group any two maximal tori are conjugate.

Let $\Delta$ in $N_\R$ be the fan of $X$, a fake weighted projective space with weights $(\lambda_0,\lambda_1,\ldots,\lambda_n)$. Let $\rho_0,\rho_1,\ldots,\rho_n$ be primitive elements of $N$ which generate the one-skeleton of $\Delta$. We have that:
\begin{equation}\label{eq:make_fake_to_wps_1}
\sum_{i=0}^n\lambda_i\rho_i=0.
\end{equation}
Let $N'\subset N$ be the lattice generated by the $\rho_i$. Let $\Delta'$ be the projection of $\Delta$ onto $N'_\R$. By construction the corresponding $\rho'_i$ of $\Delta'$ generate the lattice $N'$ and satisfy equation~\eqref{eq:make_fake_to_wps_1}. Hence, by Proposition~\ref{prop:generating_fan_unique}, $\Delta'$ is the fan of $\Proj(\lambda_0,\lambda_1,\ldots,\lambda_n)$. We obtain:
\begin{prop}\label{prop:map_fake_to_wps}
Let $X$ be any fake weighted projective space with weights $(\lambda_0,\lambda_1,\ldots,\lambda_n)$. There exists a 
finite Galois \'etale in codimension one morphism $\Proj(\lambda_0,\lambda_1,\ldots,\lambda_n)\rightarrow X$.
\end{prop}


\begin{cor}[cf.~\protect{\cite[Proposition~4.7]{Con02}}]\label{cor:fakes_are_quotients}
Let $X$ be any fake weighted projective space with weights $(\lambda_0,\lambda_1,\ldots,\lambda_n)$. Then $X$ is the quotient of $\Proj(\lambda_0,\lambda_1,\ldots,\lambda_n)$ by the action of the finite group $N/N'$ acting free in codimension one.
\end{cor}
\begin{cor}[cf.~\protect{\cite[Proposition~1.7]{Reid82}}]\label{cor:map_fake_to_wps}
Let $X$ be any fake weighted projective space with weights $(\lambda_0,\lambda_1,\ldots,\lambda_n)$. If $X$ has at worst terminal (resp. canonical) singularities then $\Proj(\lambda_0,\lambda_1,\ldots,\lambda_n)$ has at worst terminal (resp. canonical) singularities.
\end{cor}


Corollary~\ref{cor:map_fake_to_wps} tells us that if we wish to classify all fake weighted projective spaces with at worst terminal (resp. canonical) singularities, it is sufficient to find only those weights $(\lambda_0,\lambda_1,\ldots,\lambda_n)$ for which the corresponding weighted projective space possesses at worst terminal (resp. canonical) singularities. In essence, there do not exist any `extra' weights.

A similar result holds for Gorenstein fake weighted projective space:
\begin{cor}\label{cor:map_Gorenstein_to_wps}
With notation as above, $X$ is Gorenstein only if $\Proj(\lambda_0,\lambda_1,\ldots,\lambda_n)$ is Gorenstein.
\end{cor}
Let $P:=\sconv{\rho_0,\ldots,\rho_n}\subset N_\R$ be an $n$-simplex, and define the \emph{dual} by:
$$P^\vee:=\set{u\in M_\R\mid u(v)\geq -1\text{ for all }v\in P}.$$
There is a fascinating result concerning the weights of dual simplices, due to Conrads:
\begin{prop}[\protect{\cite[Lemma~5.3]{Con02}}]\label{prop:dual_weights_are_equal}
Let $X(P)$ be any Gorenstein fake weighted projective space with weights $(\lambda_0,\lambda_1,\ldots,\lambda_n)$ and associated $n$-simplex $P$. Then the fake weighted projective space $X(P^\vee)$ also has weights $(\lambda_0,\lambda_1,\ldots,\lambda_n)$.
\end{prop}
It should be noted that weights of Gorenstein weighted projective space are well understood~(see~\cite{Bat94}): A weighted projective space $\Proj(\lambda_0,\ldots,\lambda_n)$ is Gorenstein if and only if each $\lambda_j\mid\sum\lambda_i$. Hence the weights can be expressed in terms of unit partitions, and are intimately connected with the Sylvester sequence~(\cite{Nill04}).

Corollary~\ref{cor:fakes_are_quotients} provides the motivation for the following definition:
\begin{defn}\label{defn:multiplicity}
Let $P\subset N_\R$ be a $n$-simplex whose vertices $\rho_0,\rho_1,\ldots,\rho_n$ are contained in the lattice $N$. We define the \emph{multiplicity} of $P$ to be the index of the lattice generated by the $\rho_i$ in the lattice $N$; in other words, equal to the order of the group $N/N'$. We write:
$$\mathrm{mult}\,P:=[N:\Z\rho_0+\Z\rho_1+\ldots+\Z\rho_n].$$
\end{defn}
By Proposition~\ref{prop:generating_fan_unique} we have that $X(P)$ is a weighted projective space if and only if $\mathrm{mult}\, P=1$. In fact there exists a bound on how large $\mathrm{mult}\,P$ can be; this depends only on the weights and the number of interior lattice points $\abs{N\cap P^\circ}$ of $P$ (see Theorem~\ref{thm:bound_fake_wps_mult}). Before we can prove the existence of this bound, we shall require a generalisation of Minkowski's Theorem. Throughout, the volume is given relative to the underlying lattice.
\begin{thm}[\cite{Cor35}]\label{thm:Corput}
Let $k$ be any positive integer and let $K\subset N_\R$ be any centrally symmetric convex body such that $\mathrm{vol}\,K> 2^nk$.
Then $K$ contains at least $k$ pairs of points in the lattice $N$.
\end{thm}
\begin{cor}\label{cor:simplex_volume_bound}
Let $P:=\sconv{\rho_0,\rho_1,\ldots,\rho_n}\subset N_\R$ be any simplex such that:
$$\sum_{i=0}^n\lambda_i\rho_i=0,\qquad\text{for some }\lambda_i\in\Z_{>0}.$$
Let $h:=\sum_{i=0}^n\lambda_i$ and $k:=\abs{N\cap P^\circ}$. Then:
$$\mathrm{vol}\,P\leq\frac{kh^n}{n!\lambda_1\lambda_2\ldots\lambda_n}.$$
\end{cor}
\begin{proof}
Consider the convex body:
$$K:=\set{\sum_{i=1}^n\mu_i(\rho_i-\rho_0) \Big| \abs{\mu_i}\leq\frac{\lambda_i}{h}}.$$
This is centrally symmetric around the origin, with volume:
$$\mathrm{vol}\,K=\left(n!\prod_{i=1}^n\frac{2\lambda_i}{h}\right)\mathrm{vol}\,P.$$
If $\mathrm{vol}\,K>2^nk$ then, by Theorem~\ref{thm:Corput}, at least $k$ pairs of lattice points lie in the interior of $P$. But this contradicts the definition of $k$. Hence $\mathrm{vol}\,K\leq 2^nk$ and the result follows.
\end{proof}
Corollary~\ref{cor:simplex_volume_bound} can also be found in~\cite[Theorem~3.4]{Hen83},~\cite[Lemma~2.3]{LZ91}, or~\cite[Lemma~5]{Pik01}.
\begin{thm}\label{thm:bound_fake_wps_mult}
Let $P$ be the $n$-simplex associated with a fake weighted projective space $X$ with weights $(\lambda_0,\lambda_1,\ldots,\lambda_n)$. Then:
$$\mathrm{mult}\,P\leq\frac{\abs{N\cap P^\circ}h^{n-1}}{\lambda_1\lambda_2\ldots\lambda_n},\qquad\text{where }h:=\sum_{i=0}^n\lambda_i.$$
\end{thm}
\begin{proof}
Let $P'$ be the simplex associated with $\Proj(\lambda_0,\lambda_1,\ldots,\lambda_n)$, and let $F_i$ be the facet of $P'$ not containing the vertex $\rho_i'$. By considering the order of the group action on the affine patch corresponding to $F_i$, we see that $\abs{\det F_i}=\lambda_i$. Summing over all facets, we obtain:
\begin{equation}\label{eq:wps_volume}
\mathrm{vol}\,P'=\frac{h}{n!}.
\end{equation}
Combining Proposition~\ref{prop:map_fake_to_wps} with equation (\ref{eq:wps_volume}) gives:
\begin{equation}\label{eq:volume_of_fake_wps}
\mathrm{vol}\,P=\frac{h}{n!}\,\mathrm{mult}\,P.
\end{equation}
Finally we apply Corollary~\ref{cor:simplex_volume_bound}.
\end{proof}
The omission of $\lambda_0$ in the denominator is intentional. Of course it makes sense to choose the $\lambda_i$ such that $\lambda_0\leq \lambda_j$ for all $j>0$. It is reasonable to conjecture that the missing factor $\lambda_0$ should appear in the denominator, making this bound tighter.

\begin{cor}\label{cor:bound_canon_fake_wps_mult}
With notation as above, assume that $X$ has at worst canonical singularities. Then:
\begin{equation}\label{eq:canon_mult_bound}
\mathrm{mult}\,P\leq\frac{h^{n-1}}{\lambda_1\lambda_2\ldots\lambda_n}.
\end{equation}
\end{cor}
When $X$ is canonical, the right-hand side of (\ref{eq:canon_mult_bound}) is remarkably similar to the degree $(-K_\Proj)^n$ of $\Proj(\lambda_0,\lambda_1,\ldots,\lambda_n)$; the inequality becomes:
$$\mathrm{mult}\,P\leq\frac{\lambda_0}{h}(-K_\Proj)^n.$$

We conclude by mentioning two rather neat results of Conrads, for which we need the following definition.
\begin{defn}\label{defn:Hermite_normal_form}
For $n,k\in\Z_{>0}$ we denote by $\mathrm{Herm}(n,k)$ the set of all lower triangular matrices $H=(h_{ij})\in GL(n,\Q)\cap M(n\times n;\Z_{\geq 0})$ with $\det H=k$, where $h_{ij}\in\set{0,\ldots,h_{jj}-1}$ for all $j=1,\ldots,n-1$ and all $i>j$. We call $\mathrm{Herm}(n,k)$ then set of \emph{Hermite normal forms} of dimension $n$ and determinant $k$.
\end{defn}
\begin{thm}[\protect{\cite[Theorem~4.4]{Con02}}]\label{thm:equiv_Hermite}
Let $X(P')$ be any fake weighted projective space with weights $(\lambda_0,\lambda_1,\ldots,\lambda_n)$ and associated $n$-simplex $P'$. Let $P$ the $n$-simplex associated with $\Proj(\lambda_0,\lambda_1,\ldots,\lambda_n)$. Then there exists $H\in\mathrm{Herm}(n,\mathrm{mult}\,P')$ such that $P'=HP$ (up to the action of $GL(n,\Z)$).
\end{thm}
\begin{cor}[\protect{\cite[Proposition~5.5]{Con02}}]\label{cor:reflexive_multiplicities_divide}
With notation as above, if $X(P')$ is Gorenstein then:
$$\mathrm{mult}\,P'\mid\mathrm{mult}\,P^\vee.$$
\end{cor}
\begin{proof}
Since $P'$ is reflexive, so $P$ must be reflexive by Corollary~\ref{cor:map_Gorenstein_to_wps}. By Theorem~\ref{thm:equiv_Hermite} there exists some $H\in\mathrm{Herm}(n,\mathrm{mult}\,P')$ such that $P'=HP$. Hence $P'^\vee=H^\vee P^\vee$. Now $H^\vee=(H^t)^{-1}$, and so $\det H^\vee=1/\mathrm{mult}\,P'$.

Thus $\det P'^\vee=\det P^\vee/\mathrm{mult}\,P'$. By Proposition~\ref{prop:dual_weights_are_equal} and equation~\eqref{eq:volume_of_fake_wps} we obtain:
$$\mathrm{mult}\,P'^\vee=\frac{\mathrm{mult}\,P^\vee}{\mathrm{mult}\,P'}.$$
Observing that $\mathrm{mult}\,P'^\vee\in\Z_{>0}$ gives the result.
\end{proof}
\section{Upper Bounds on the Weights}\label{sec:upper_bound_on_bary}
Let $X$ be a fake weighted projective space with weights $(\lambda_0,\ldots,\lambda_n)$, where $\lambda_0\leq\ldots\leq\lambda_n$ and $\gcd\set{\lambda_0,\ldots,\lambda_n}=1$. Throughout we shall assume that $X$ has at worst canonical singularities.

Let $P:=\sconv{\rho_0,\ldots,\rho_n}\subset N_\R$ be the associated simplex. By assumption there is a unique interior lattice point of $P$, namely the origin, and:
$$\sum_{i=0}^n\frac{\lambda_i}{h}\rho_i=0,\qquad\text{ where }h:=\sum_{i=0}^n\lambda_i.$$

In~\cite{Pik01} an upper bound is given for the volume of $P$:
\begin{thm}[\protect{\cite[Theorem~6]{Pik01}}]\label{thm:bound_Fano_simplex}
With notation as above, we have:
$$\mathrm{vol}\,P\leq\frac{1}{n!}2^{3n-2}15^{(n-1)2^{n+1}}.$$
\end{thm}
Combining this result with equation (\ref{eq:wps_volume}) immediately gives us an upper bound on $h$. Unfortunately this bound is far from tight. In the case when $X$ is Gorenstein, combining equation (\ref{eq:wps_volume}) with~\cite[Theorem~C]{Nill04} provides much better bounds:
\begin{prop}\label{prop:Nills_bound_on_h}
Suppose that $X$ is Gorenstein. Then $h\leq t_n$, where $t_n:=y_n-1$ is defined in terms of the Sylvester sequence $y_0:=2$, $y_k:=1+y_0\cdots y_{k-1}$.
\end{prop}
A lower bound on $\lambda_0/h$ was also presented in~\cite{Pik01}:
\begin{thm}[\protect{\cite[Theorem~2]{Pik01}}]\label{thm:lambda_0_lower_bound}
With notation as above;
$$\frac{\lambda_0}{h}\geq\frac{1}{8\cdot 15^{2^{n+1}}}.$$
\end{thm}
When $X$ is Gorenstein, \cite[Proposition~3.4]{Nill04} establishes the following lower bounds:
\begin{prop}
Suppose that $X$ is Gorenstein. With notation as above, for any $k\in\set{0,\ldots,n}$ we have that:
$$\frac{\lambda_{k}}{h}\geq\frac{1}{(k+1)t_{n-k}}.$$
\end{prop}
We shall prove the following upper bounds hold:
\begin{thm}\label{thm:barycentric_bounds}
With notation as above, for any $k\in\set{2,\ldots,n}$ we have that:
$$\frac{\lambda_{k}}{h}\leq\frac{1}{n-k+2},$$
with strict inequality if $X$ possesses at worst terminal singularities.
\end{thm}
We shall require the following elementary lemma:
\begin{lem}\label{lem:negative_cone_and_origin}
Let $\sigma=\scone{x_1,\ldots,x_m}$ be an $m$-dimensional convex cone. If $x\in-\sigma$ then $0\in\sconv{x,x_1,\ldots,x_m}$.
\end{lem}
We shall present our proof of Theorem~\ref{thm:barycentric_bounds} assuming that $X$ is terminal; the result when $X$ is canonical should be apparent.
\begin{proof}[Proof of Theorem~\ref{thm:barycentric_bounds}]
Since $\sum_{i=0}^n\lambda_i\rho_i=0$, so:
$$\sum_{i=0}^{n-k-1}\lambda_i\rho_i=\sum_{j=n-k}^n-\lambda_j\rho_j.$$
Now $\sum_{i=0}^{n-k-1}\lambda_i=h-\sum_{j=n-k}^n\lambda_j$, giving:
$$x:=\sum_{j=n-k}^n\frac{-\lambda_j}{h-l}\rho_j\in\sconv{\rho_0,\ldots,\rho_{n-k-1}},\qquad\text{where }l:=\sum_{j=n-k}^n\lambda_j.$$
Since $P$ is simplicial, $\sconv{\rho_0,\ldots,\rho_{n-k-1}}$ is a face of $P$. Since the $\lambda_i$ are all strictly positive, $x$ lies strictly in the interior of this face.

Let us suppose for a contradiction that:
\begin{equation}\label{eq:barycentric_bounds_1}
\lambda_{n-k+i}\geq\frac{h}{k+2},\qquad\text{for all }i\in\set{0,\ldots,k}.
\end{equation}
Consider the $(k+1)$-dimensional lattice $\Gamma$ generated by $e_0,\ldots,e_k$. There exists a map of lattices $\gamma:\Gamma\rightarrow N$ given by sending $e_i\mapsto \rho_{n-k+i}$. Note that this map is injective. Let $x':=\sum_{i=0}^k-\lambda_{n-k+i}/(h-l)e_i$. We shall show that the non-zero lattice point $p:=-\sum_{i=0}^ke_i$ lies in $\sconv{x',e_0,\ldots,e_k}$. Hence $\gamma(p)\neq 0$ is a lattice point in $\sconv{x,\rho_{n-k},\ldots,\rho_n}\subset P$.

Since $p\notin\sconv{e_0,\ldots,e_k}$, so $\gamma(p)$ is not contained in $\sconv{\rho_{n-k},\ldots,\rho_n}$. The only remaining possibility which does not contradict $P$ having only the origin as a strictly internal lattice point is that $\gamma(p)=x$. But if $P$ is terminal we have a contradiction.

Consider $\lambda_n$. By~\eqref{eq:barycentric_bounds_1} we have that:
\begin{equation}\label{eq:barycentric_bounds_1a}
\lambda_n-h\geq\frac{-h(k+1)}{k+2}.
\end{equation}
Summing~\eqref{eq:barycentric_bounds_1} over $0\leq i<k$ gives:
\begin{equation}\label{eq:barycentric_bounds_1b}
l-\lambda_n\geq\frac{hk}{k+2}.
\end{equation}
Combining equations~\eqref{eq:barycentric_bounds_1a} and~\eqref{eq:barycentric_bounds_1b} gives us that $l-h\geq -h/(k+2)$. Observing that $l-h<0$, we obtain $(k+2)/h\leq 1/(h-l)$. Thus, for any $j\in\set{n-k,\ldots,n}$, we have that:
$$-1\geq\frac{-\lambda_j}{h-l},$$
i.e. the coefficients of $x'$ are all $\leq -1$.

Let $\tau$ be the lattice translation of $\Gamma$ which sends $p$ to $0$. Applying Lemma~\ref{lem:negative_cone_and_origin} to $\scone{\tau e_0,\ldots,\tau e_k}$, if $\tau(x')\in-\scone{\tau e_0,\ldots,\tau e_k}$ then $p\in\sconv{x',e_0,\ldots,e_k}$ and we are done. Hence assume that this is not the case.

Let $H_i\subset\Gamma_\R$ be the hyperplane containing the $k+1$ points $e_0,\ldots,\widehat{e_i},\ldots,e_k,$ and $p$; let $H_i^+$ be the half--space in $\Gamma_\R$ whose boundary is $H_i$ and which contains the point $2p$. Then:
$$-\scone{\tau e_0,\ldots,\tau e_k}=\tau\left(\bigcap_{i=0}^kH_i^+\right).$$
Since $\tau(x')\notin-\scone{\tau e_0,\ldots,\tau e_k}$, we have that $x'\notin H_i^+$ for some $i$. Assume, with possible reordering of the indices, that $x'\notin H_0^+$.

$H_0$ is given by:
$$\set{\sum_{j=1}^k\mu_je_j-\left(1-\sum_{j=1}^k\mu_j\right)\sum_{i=0}^ke_i\ \Big|\ \mu_i\in\R}.$$
Let $q:=\sum_{i=0}^k\nu_ie_i$ be any point in $\Gamma_\R$. By projecting $q$ onto $H_0$ along $e_0$ (regarded as a vector) we can always choose our $\mu_i$ such that:
\begin{equation}\label{eq:barycentric_bounds_2a}
\mu_j+\sum_{i=1}^k\mu_i-1=\nu_j,\qquad\text{for }1\leq j\leq k.
\end{equation}
Comparing the sign of $\sum_{i=1}^k\mu_i-1$ with $\nu_0$ tells us on which side of the hyperplane $H_0$ the point $q$ lies.

We have that $2p$ lies on the opposite side of $H_0$ to $x'$. Setting $\nu_j=-2$ for all $j$ in equation~\eqref{eq:barycentric_bounds_2a} tells us that:
$$\sum_{i=1}^k\mu_i=\frac{-k}{k+1}.$$
Hence we see that:
$$\sum_{i=1}^k\mu_i-1=\frac{-k}{k+1}-1>-2.$$
We thus require that:
\begin{equation}\label{eq:barycentric_bounds_3}
\sum_{i=1}^k\mu_i-1<\frac{-\lambda_{n-k}}{h-l}.
\end{equation}
(In other words $x'$ lies on the opposite side of $H_0$ to $2p$.)

Comparing coefficients with $x'$, we see that:
\begin{equation}\label{eq:barycentric_bounds_2}
\mu_j+\sum_{i=1}^k\mu_i-1=\frac{-\lambda_{n-k+j}}{h-l},\qquad\text{for }1\leq j\leq k.
\end{equation}
Summing equation~\eqref{eq:barycentric_bounds_2} for all $1\leq j\leq k$ and combining this with~\eqref{eq:barycentric_bounds_3} gives:
$$\sum_{j=1}^k\frac{-\lambda_{n-k+j}}{h-l}+k<\frac{-(k+1)\lambda_{n-k}}{h-l}+k+1.$$
Simplifying, and recalling that $\sum_{j=n-k}^n\lambda_j=l$, gives us that:
\begin{equation}\label{eq:barycentric_bounds_4}
\lambda_{n-k}<\frac{h}{k+2}.
\end{equation}
Equation~\eqref{eq:barycentric_bounds_4} contradicts~\eqref{eq:barycentric_bounds_1}, concluding the proof.
\end{proof}
\bibliographystyle{amsalpha}

\begin{thebibliography}{BCF{\etalchar{+}}05}
\bibitem[Bat94]{Bat94}
Victor~V. Batyrev, \emph{Dual polyhedra and mirror symmetry for {C}alabi-{Y}au
  hypersurfaces in toric varieties}, J. Algebraic Geom. \textbf{3} (1994),
  no.~3, 493--535.

\bibitem[BB92]{BB92}
A.~A. Borisov and L.~A. Borisov, \emph{Singular toric {F}ano three-folds}, Mat.
  Sb. \textbf{183} (1992), no.~2, 134--141, text in Russian. English transl.:
  \emph{Russian Acad. Sci. Sb. Math.}, \textbf{75} (1993), 277--283.

\bibitem[BCF{\etalchar{+}}05]{Snow05}
Matthias Beck, Beifang Chen, Lenny Fukshansky, Christian Haase, Allen Knutson,
  Bruce Reznick, Sinai Robins, and Achill Sch{\"u}rmann, \emph{Problems from
  the {C}ottonwood {R}oom}, Integer points in polyhedra---geometry, number
  theory, algebra, optimization, Contemp. Math., vol. 374, Amer. Math. Soc.,
  Providence, RI, 2005, pp.~179--191.

\bibitem[Buc02]{Krych02}
Weronika Buczy{\'n}ska, \emph{Toryczne przestrzenie rzutowe}, Magister thesis (2002),
  text in Polish, available from \href{http://www.mimuw.edu.pl/~jarekw/}
  {\texttt{http://www.mimuw.edu.pl/$\sim$jarekw/}}. English
  transl.: \href{http://arxiv.org/abs/0805.1211} {\texttt{arXiv:0805.1211}} (2008).

\bibitem[Con02]{Con02}
Heinke Conrads, \emph{Weighted projective spaces and reflexive simplices},
  Manuscripta Math. \textbf{107} (2002), no.~2, 215--227.

\bibitem[Dem70]{Dem70}
Michel Demazure, \emph{Sous-groupes alg\'ebriques de rang maximum du groupe de
  {C}remona}, Ann. Sci. \'Ecole Norm. Sup. (4) \textbf{3} (1970), 507--588.

\bibitem[Fuj03]{Fuj03}
Osamu Fujino, \emph{Notes on toric varieties from {M}ori theoretic viewpoint},
  Tohoku Math. J. (2) \textbf{55} (2003), no.~4, 551--564.

\bibitem[Hen83]{Hen83}
Douglas Hensley, \emph{Lattice vertex polytopes with interior lattice points},
  Pacific J. Math. \textbf{105} (1983), no.~1, 183--191.

\bibitem[Kas06a]{Kas03}
Alexander~M. Kasprzyk, \emph{Toric {F}ano threefolds with terminal
  singularities}, Tohoku Math. J. (2) \textbf{58} (2006), no.~1, 101--121.

\bibitem[Kas06b]{KasPhD}
\bysame, \emph{Toric {Fano} varieties and convex polytopes}, PhD thesis,
  University of Bath (2006), available from
  \href{http://hdl.handle.net/10247/458}{\texttt{http://hdl.handle.net/10247/4%
58}}.

\bibitem[LZ91]{LZ91}
Jeffrey~C. Lagarias and G{\"u}nter~M. Ziegler, \emph{Bounds for lattice
  polytopes containing a fixed number of interior points in a sublattice},
  Canad. J. Math. \textbf{43} (1991), no.~5, 1022--1035.

\bibitem[Mat02]{Mat02}
Kenji Matsuki, \emph{Introduction to the {M}ori program}, Universitext,
  Springer-Verlag, New York, 2002.

\bibitem[Nil07]{Nill04}
Benjamin Nill, \emph{Volume and lattice points of reflexive simplices},
  Discrete Comput. Geom. \textbf{37} (2007), no.~2, 301--320.

\bibitem[Pik01]{Pik01}
Oleg Pikhurko, \emph{Lattice points in lattice polytopes}, Mathematika
  \textbf{48} (2001), no.~1-2, 15--24 (2003).

\bibitem[Rei80]{Reid82}
Miles Reid, \emph{Canonical {$3$}-folds}, Journ\'ees de G\'eometrie
  Alg\'ebrique d'Angers, Juillet 1979/Algebraic Geometry, Angers, 1979,
  Sijthoff \& Noordhoff, Alphen aan den Rijn, 1980, pp.~273--310.

\bibitem[Rei83]{Reid83T}
\bysame, \emph{Decomposition of toric morphisms}, Arithmetic and geometry, Vol.
  II, Progr. Math., vol.~36, Birkh\"auser Boston, Boston, MA, 1983,
  pp.~395--418.

\bibitem[Rei87]{Reid85}
\bysame, \emph{Young person's guide to canonical singularities}, Algebraic
  geometry, Bowdoin, 1985 (Brunswick, Maine, 1985), Proc. Sympos. Pure Math.,
  vol.~46, Amer. Math. Soc., Providence, RI, 1987, pp.~345--414.

\bibitem[vdC35]{Cor35}
J.~G. van~der Corput, \emph{Verallgemeinerung einer {M}ordellschen
  {B}eweismethode in der {G}eometrie der {Z}ahlen}, Acta Arithm. \textbf{1}
  (1935), 62--66, text in German.
\end{thebibliography}
\newcommand{\etalchar}[1]{$^{#1}$}
\providecommand{\bysame}{\leavevmode\hbox to3em{\hrulefill}\thinspace}
\providecommand{\MR}{\relax\ifhmode\unskip\space\fi MR }
\providecommand{\MRhref}[2]{%
  \href{http://www.ams.org/mathscinet-getitem?mr=#1}{#2}
}
\providecommand{\href}[2]{#2}

\end{document}